\documentclass[11pt]{article}

\usepackage{amssymb}
\usepackage{amsfonts}
\usepackage{amsmath}
\usepackage[english,german]{babel}
\usepackage{epsfig}
\usepackage{latexsym}
\usepackage{color}
\usepackage{epsfig}
\usepackage{psfrag,subfigure}

\newtheorem{thm}{Theorem}[section]
\newtheorem{prop}[thm]{Proposition}
\newtheorem{cor}[thm]{Corollary}
\newtheorem{lemma}[thm]{Lemma}
\newtheorem{definition}[thm]{Definition}
\newtheorem{rem}[thm]{Remark}

\newtheorem{Assumption}[thm]{Assumption}

\newenvironment{Proof}{\textsc{Proof.}}{\mbox{ } \hfill $\Box$ \vspace{2mm}}

\newcommand{\ba}{\begin{array}{ll}}
\newcommand{\bal}{\begin{array}{ll}}
\newcommand{\ea}{\end{array}}

\newcommand{\n}{\mathbb{N}}
\newcommand{\re}{\mathbb{R}}
\newcommand{\E}{\mathbb{E}}
\newcommand{\Prob}{\mathbb{P}}

\def\dfrac{\displaystyle\frac}

\def\n{\mbox{$\mathbb{N}$}}

\def\X{\mbox{$\cal{X}$}}
\def\t{\mbox{${\theta}$}}
\def\x{\mbox{${\theta}$}}
\def\T{\mbox{${\Theta}$}}
\def\O{\mbox{${\Omega}$}}
\def\u{\overline{u}}
\def\f{\overline{f}}
\def\v{\overline{v}}

\def\D{\overline{D}}

\def\C{\cal{C}}

\begin{document}

\selectlanguage{english}

\renewcommand{\baselinestretch}{1.25} \normalsize


\title{A Numerical Approach to the Estimation of Solutions of some Variational Problems with Convexity Constraints \thanks{We thank Guillaume Carlier and Yves
Lucet for their thoughtful comments and suggestions.}}

\author{\normalsize Ivar Ekeland \\[8pt]
        \small  Department of Mathematics  \\
        \small University of British Columbia \\
        \small 1984 Mathematics Road\\
        \small Vancouver, BC, V6T 1Z2\\
        \small ekeland@math.ubc.ca
         \and
         \normalsize  Santiago Moreno-Bromberg \\[8pt]
        \small Department of Mathematics\\
        \small University of British Columbia \\
        \small 1984 Mathematics Road\\
        \small Vancouver, BC, V6T 1Z2\\
        \small smoreno@math.ubc.ca
\vspace*{0.8cm}}

\maketitle

\begin{abstract}
We present an algorithm to approximate the solutions to variational problems where set of admissible functions consists of convex functions. The main motivator behind this numerical method is estimating  solutions to Adverse Selection problems within a Principal-Agent framework. Problems such as product lines design, optimal taxation, structured derivatives design, etc. can be studied through the scope of these models. We develop a method to estimate their optimal pricing schedules.
\end{abstract}

\vspace{1mm}

\begin{center}
\textsc{Preliminary - Comments Welcome}
\end{center}

\vspace{1mm}


\vspace{2mm}


\hfill

\newpage

\renewcommand{\baselinestretch}{1.20} \normalsize

\section{Introduction}
Arguably, Newton's problem of the body of minimal resistance is the original variational problem with convexity constraints. It consists of finding the shape of a solid that encounters the least resistance when moving through a fluid. This is equivalent to finding a convex function from a convex domain (originally a disk)  in $\re^2$ to $\re$ that minimizes a certain funtional (see section \ref{problem}). Newton's original ``solution" assumed radial symmetry. This turned out to be false, as shown by Brock, Ferone and Kawohl in \cite{kn:bfk}, which sparked new interest to the study of variational problems with convexity constraints. One can also find these kinds of problems in finance and economics. Starting in 1978 with the seminal paper of Mussa and Rosen \cite{kn:mr}, the study of non-linear pricing as a means of market screening under Adverse-Selection has produced a considerable stream of contributions (\cite{kn:A},\cite{kn:cet},\cite{kn:rch},...). In models where goods are described by a single {\rm quality} and the set of {\rm agents} is differentiated by a single parameter, it is in general possible to find closed form solutions for the {\rm pricing schedule}. This is, however, not the case when multidimensional {\rm consumption bundles} and agent {\rm types} are considered. Although Rochet and Chon\'e \cite{kn:rch} provided conditions for the existence of an optimal pricing rule  and fully characterized the ways in which markets differentiate in a multidimensional setting, they  also pointed out that it is only in very special cases that one can expect to find closed form solutions.  The same holds true for models where the set of goods lies in an infinite-dimensional space, even when agent types are one-dimensional. This  framework was first used, to our knowledge, by  Carlier, Ekeland and Touzi \cite{kn:cet} to price financial derivatives traded ``over-the-counter". It was then extended by Horst and Moreno \cite{kn:HM} to model the actions of a monopolist who has an initial risky position that she evaluates via a coherent risk measure, and who intends to transfer part of her risk to a set of heterogenous agents. In both cases the authors find that only very restrictive examples allow for explicit solutions.

Given that a great variety of problems, such as product lines design, optimal taxation, structured derivatives design, etc. can be studied through
the scope of these models, there is a clear need for robust and efficient numerical methods that approximate their optimal pricing schedules. Note that this also provides an approximation of the optimal ``products". Most of the papers mentioned above eventually face solving a variational problem under convex constraints. This family of problems have lately been studied under different scopes. Carlier and Lachand-Robert \cite{kn:CLR} have studied the $C^1$ regularity of minimizers when the functional is elliptic and the admissible functions satisfy a Dirichlet-type boundary condition. Their results can be extended to of our examples. Lachand-Robert and Pelletier \cite{kn:LRP} characterize the extreme points of a functional depending only on $\nabla f$ over a set of convex functions with uniform convex bounds. In this paper we provide several variants of an algorithm, based on the idea of approximating a convex function by an affine envelope, to solve these types of problems. This deviates from previous work by Chon\'e  and Herv\'e \cite{kn:ch} and Carlier, Lachand-Robert and Maury \cite{kn:cl}, where the authors use finite element methods. In the former case, a conformal (interior) method is used and a non-convergence result is given. As a consequence, the latter uses an exterior approximation method, which is indeed found to be convergent in the classical projection problem in $H_0^1.$ Lachand-Robert and Oudet present in  \cite{kn:LRO} and algorithm for minimizing functionals within convex bodies that shares some similarities to ours. For a particular problem, they start with an admissible polytope and iteratively modify the normals to the facets in order to find an approximate minimizer.


We estimate the minimizers for several problems with known, closed form solutions as a means of comparing the output of our method to the true solutions. These are taken from \cite{kn:rch} and \cite{kn:cet}. Finally, we provide an example in which we approximate the solution to a risk-minimization problem similar to the one presented in \cite{kn:HM}. This is still based on the affine-envelope idea, but requires some additional methodology, since it involves solving  a non-standard variational problem.

The remainder of this paper is organized as follows. In Section
\ref{problem} we state our problem and provide some classical examples.
Our algorithm and a proof of its convergence are presented in Section \ref{DescriptionConvex}. In Section \ref{Examples} we show the solutions
obtained via our algorithm to several problems found in the literature. Since these problems share a common microeconomic motivation, we include a brief discussion on the latter. The examples include the well known ``Rochet-Chon\'e" problem, a one dimensional example from Carlier, Ekeland and Touzi and  the risk transfer case for a principal who offers call options with type-dependent strikes and evaluates her risk via the  ``short fall" of her position. This section is followed by our conclusions. Finally a section devoted to technical results and all our codes are included in the appendix.

\section{Setting}\label{problem}
The aim of this paper is to present a numerical algorithm to approximate the solutions of some variational problems subject to convexity constraints. A classical example of the latter is Newton's problem of the body of minimal resistance, which, given $\T$ a smooth subset of $\re^2,$ consists of minimizing
\[
I[v]=\int_{\T}\frac{d\t}{1+|\nabla v|^2},
\]
over the set of convex functions $\{f:\T\to\re\}.$  We use the following notation throughout:
\begin{itemize}
\item $\T,Q\subset\re^n$ are convex and compact sets,

\item $L(\t,z,p)=z-\t\cdot p+C(p),$ where $C$ is strictly convex and $C^1.$

\item $\C:=\{f:\T\to\re\,\mid\,f\ge 0\,\,{\mbox{is convex, and}}\,\,\nabla f\in Q\,\,{\mbox{a.e}}\},$

\item $I[f]:=\int_{\Theta}L(\t, f(\t), \nabla f(\t))d\t.$
\end{itemize}
Our objective is to (numerically) estimate the solution to
\begin{equation}\label{eq:problem}
{\cal{P}}:=\inf_{f\in C}I[f]
\end{equation}
We assume $C$ is such that (\ref{eq:problem}) has a unique solution (see, for example, \cite{kn:GH}). Given the properties of $L$ and ${\cal{C}}$ we immediately have the following

\begin{prop}\label{zerov} Assume $\v$ solves ${\cal{P}},$ then there is $\t_0$ in $\T$ such that $\v(\t_0)=0.$
\end{prop}
\begin{Proof}
Let $\v_0=\min_{\t\in\T}\v(\t)$ (recall $\T$ is compact)and define $\overline{u}(\t):=\v(\t)-\v_0,$ then
\[
I[\u]=\int_{\T}\u(\t)-\t\cdot\nabla\u(\t)+C(\nabla\u(\t))d\t=I[\v]-\|\T\|\v_0.
\]
This would contradict the hypothesis of $\v$ being a minimizer of $I$ over $\C$ unless $\v_0=0.$
\end{Proof}

It follows from proposition \ref{zerov} that we can redefine $\C$ to include only functions that have a root in $\T.$ This, together with the compactness of $Q,$ implies the following proposition, which we will use frequently.

\begin{prop}\label{zerov2} There exists $0<K<\infty$ such that $v\le K$ for all $v$ in $\C.$
\end{prop}
It follows from Proposition \ref{zerov2} and the restriction on the gradients that for each choice of function $C,$ problem ${\cal{P}}$ has a unique solution, since the functional $I$ will be strictly convex, lower semi continuous and the admissible set is bounded (see \cite{kn:et}).

\begin{rem}
Our algorithm will still work with more general $L$'s as long as one can prove that the family of feasible minimizers is uniformly bounded.
\end{rem}

\section{Description of the Algorithm}\label{DescriptionConvex}
From this point on, whenever we use supscripts we refer to vectors. For example $V^k=(V^k_1,\ldots, V^k_m).$ On the other hand a subscript indicates a function to be evaluated over some closed, convex subset of $\re^k$ of non-empty interior, ie, $\{V_k\}$ is a sequence of functions $V_k:X\to\re$ for some $X$ contained in $\re^n.$

\begin{Assumption} We will consider $\T=[a,b]^n.$
\end{Assumption}

To find an approximate solution to ${\cal{P}},$ we proceed as follows:

\begin{enumerate}

\item We discretize the domain $\T$ in the following way: We partition it into $\Sigma_k,$ which consists of $k^n$ equal cubes of volume $\|\Sigma_k\|:=\left(\frac{b-a}{k}\right)^n.$ The elements of $\Sigma_k$ will be denoted by $\sigma_j^k,$ $1\le j\le k^n.$ Now define $\T_k$ as the set of centers of the $\sigma_j^k$'s. The elements of $\T_k$ will be denoted by $\t_j^k.$ The choice of a uniform partition is done for computational simplicity.

\item We denote $f_i=f\left(\t_i^k\right)$ and associate such weight with $\t_i^k.$

\item  We associate to each element $\t_i^k$ of $\T_k$  a non-negative number $v_i^k$ and
an n-dimensional vector $D_i^k.$ The former represents the value of $v(\t_i^k)$ and the latter $\nabla v(\t_i^k).$

\item We solve the (non-linear) program

\begin{equation}\label{eq:ppd}
{\cal{P}}_k:=\inf\|\Sigma_k\|\sum_{i=1}^{k^n}L\left(\t_i, v_i, D_i\right)f_i
\end{equation}

\noindent over the set of all vectors of the form  $v=(v_1,\ldots,v_{k^n})$ and all matrices of the form
$D=(D_1,\ldots,D_{k^n})$ such that:

\begin{enumerate}
\item $v\ge 0$ (non-negativity),

\item $D_i\in Q$ for $i=1,\ldots k^n$ (feasibility),

\item $v_i-v_j+D_i\cdot(\t_j-\t_i)\le 0$ (convexity).
\end{enumerate}

If the problem in hand includes Dirichlet boundary conditions these can be included here as linear constraints that the $D_i$'s corresponding to points on the ``boundary" of $\T_k$ must satisfy.

\item Let $(\v^k, \D^k)$ be the solution to ${\cal{P}}_k.$ We define $\v_k(\t):=\max_i p_i(\t),$ where

\[
p_i(\t)=\v^k_i+\D^k_i\cdot (\t-\t_i).
\]

\item $\v_k$ yields an approximation to the minimizer of ${\cal{P}}.$

\end{enumerate}

\begin{rem}\label{rem:constraints} The constraints of the non-linear program determine a convex set.
\end{rem}

\begin{rem}\label{rem:pwafine} 4 (c) guarantees that $p_i$ is a supporting hyperplane of
the convex hull of the points $\{(\t_1, v_1),\ldots,(\t_{k^n}, v_{k^n})\}.$ Note that $\v_k$ is a piecewise affine convex function.
\end{rem}

\subsection{Convergence of the Algorithm}\label{ConvergenceConvex}

\begin{prop} Under the assumptions made on $L,$ the problem ${\cal{P}}_k$ has a unique solution.\end{prop}
\begin{Proof} The function
\[
J_k(v^k,D^k):=\sum_{i=1}^{k^n}\left(\t_i\cdot D_i^k-v_i^k-C(D_i^k)\right)\|\Sigma_k\|f_i
\]
is strictly convex. It follows from proposition \ref{zerov2} that any acceptable vector-matrix pair $(v^k,D^k)$
must lie in $[0,K]^k\times Q^k,$ which together with Remark
\ref{rem:constraints} implies ${\cal{P}}_k$ consists of minimizing a
strictly convex function over a compact and convex set. The result then
follows from general theory.
\end{Proof}

\begin{prop}\label{pr:conv}  There exists $\v\in{\cal{C}}$ such that:
\begin{enumerate}
\item The sequence $\{\v_k\}$ generated by the ${\cal{P}}_k$'s has a
subsequence  $\{\v_{k_j}\}$that converges uniformly to
$\overline{v}.$

\item $\lim_{k_j\to\infty} I[\v_{k_j}]=I[\overline{v}].$
\end{enumerate}
\end{prop}
\begin{Proof} The bounded (Proposition \ref{zerov2}) family $\{\v_j\}$ is  uniformly equicontinuous, as it consists of
convex functions with uniformly bounded subgradients. By the
Arzela-Ascoli theorem we have that, passing to a subsequence if
necessary, there is a non-negative and convex function
$\overline{v}$ such that

\[
\v_k\to\v\quad{\mbox{uniformly on}}\,\,\T.
\]
By convexity $\nabla\v_k\to\nabla\v$ almost everywhere (lemma \ref{pr:vprime}); since $\nabla\v_k(\t)$ belongs to the bounded set $Q,$ the integrands are dominated. Therefore, by Lebesgue Dominated Convergence we have

\[
\lim_{k\to\infty}I\left[\v_k\right]=I[\v] .
\]
\end{Proof}

Let $\u$ be the maximizer of $I[\cdot]$ within $C.$  Our aim is to
show that $\{\v_k\}$ is a minimizing sequence of problem ${\cal{P}},$ in
other words that

\[
\lim_{k\to\infty}I\left[\v_k\right]=I\left[\u\right] .
\]
We need the following

\begin{definition} Let $\u$ be such that $\inf_{u\in{\cal{C}}}I[u]=I[\u].$ Given the lattice $\T_k,$ we define:
\begin{enumerate}
\item  $\u^k_i:=\overline{u}(\t_i),$

\item $G^k_i:=\nabla {u}(\t_i),$

\item $q_i(\t):=\u^k_i+G^k_i\cdot(\t-\t_i)$ and

\item $\u_k(\t):=\sup_i q_i(\t).$
\end{enumerate}
\end{definition}

Notice that $\u_k(\t)$ is also constructed as the convex envelope of a family of affine functions. The inequalities

\begin{equation}\label{eq:JN}
 J_k(\u^k, G^k)\ge J_k(\v^k, \D^k)
\end{equation}

\begin{equation}\label{eq:IN}
I[\v_k]\ge I[\u]
\end{equation}

\noindent follow from the definitions of  $J_k(\v^k, \D^k),$ $\v_k$ and
$\u_k,$ as does the following

\begin{prop}\label{pr:conv} Let $\u$ and $\u_k$ be as above, then $\u_k\to\u$ uniformly as $k\to\infty.$\end{prop}

\begin{prop}\label{pr:conv2} For each $k$ there exist $\epsilon_1(k)$ and $\epsilon_2(k)$  such that

\begin{equation}\label{eq:JI1}
\left| J_k(\v^k, \D^k)-I[\v_k]\right|\le\epsilon_1(k)
\end{equation}

\begin{equation}\label{eq:JI2}
\left| J_k(\u^k, G^k)-I[\u_k]\right|\le\epsilon_2(k)
\end{equation}

\noindent and $\epsilon_1(k),\epsilon_2(k)\to 0$ as $k\to\infty.$
\end{prop}
\begin{Proof}
We will show (\ref{eq:JI1}) holds, the proof for (\ref{eq:JI2}) is
analogous. Define the simple function
\[
w_k(\t):=L(\t_j^k, v_j^k, D_j^k),\,\, \t\in \sigma_j^k,
\]
hence

\begin{equation}J_k(\v^k, \D^k)=\int_{\T}w_k(\t) d\t. \end{equation}
The left-hand side of (\ref{eq:JI1}) can be written as

\begin{equation}\label{eq:JI3}
\left|\int_{\T}w_k(\t)d\t-I[\v_k]\right|
\end{equation}
It follows from Lemma \ref{eq:affineconv} that there exists $\epsilon_1(k),$  such that
\[
\left|\int_{\T}w_k(\t)d\t-I[v_k]\right|\le\epsilon_1(k)
\]
and
\[
\epsilon_1(k)\to 0\quad k\to\infty.
\]

\end{Proof}

\noindent We can now prove the main theorem in this section, namely

\begin{thm} The sequence $\{\v_k\}$ is minimizing for problem ${\cal{P}}.$\end{thm}
\begin{Proof}
It follows from Proposition \ref{pr:conv2} and equation (\ref{eq:JN}) that

\begin{equation}\label{eq:conv3}
I[\u_k]+\epsilon_2(k)+\epsilon_1(k)\ge I[\v_k]\ge I[\u]
\end{equation}
Letting $k\to\infty$ in (\ref{eq:conv3}) and using
Proposition \ref{pr:conv} yields the desired result.
\end{Proof}
\section{Examples}\label{Examples}
In this section we show some results of implementing our algorithm.
The first two examples reduce quadratic programs, whereas the
third and fourth ones  are non-linear optimization programs. 
All the computer coding has been written in MatLab. However, in both cases supplemental Optimization Toolboxes were used. In the first two examples we used the Mosek 4.0 Optimization Toolbox, wherease in the last two we used Tomlab 6.0.
These four examples share a common microeconomic motivation, for which we provide an overview. We refer the interested reader to \cite{kn:AS} for a comprehensive presentation of Principal-Agent models and Adverse Selection, as well as multiple references.

\subsection{Some Microeconomic Motivation}\label{Setup}

Consider an economy with a single {\sc principal} and a continuum
of {\sc agents}. The latter's preferences are characterized by n-dimensional
vectors. These are called the agents'
{\sc types}. The set of all types will be denoted by
$\T\subset\re^n.$ The individual types $\t$ are private information,
but the principal knows their statistical distribution, which has a (non-atomic)
density $f(\t).$

Our model takes a {\sc hedonic} approach to product differentiation. We
assume goods are characterized by (n-dimensional) vectors describing
their utility-bearing attributes. The set of {\sc technologically
feasible} goods that the principal can deliver will be denoted by
$Q\subset\re_+^n,$ and it will be assumed to be compact and convex.
The cost to the principal of producing one unit
of product $p$ is denoted by $C(p).$ Products are offered on a
take-it-or-leave-it basis, each agent can buy one or zero units of a
single product $p$ and it is assumed there is no second-hand market.
The (type-dependent) preferences of the agents are represented by
the function

\[
U:\T\times Q\to\re.
\]
The (non-linear) price schedule for the technologically feasible
goods is represented by

\[
\pi:Q\to\re.
\]
When purchasing good $q$ at a price $\pi(q)$ an agent of
type $\t$ has net utility

\[
U(\t, q)-\pi(q)
\]
Each agent solves the problem

\[
\max_{q\in Q}\left\{ U(\t, q)-\pi(q)\right\}.
\]
By analyzing the choice of each agent type under a given
price schedule $\pi,$ the principal screens the market. Let

\begin{equation}\label{Max-ut}
v(\t):=U(\t, q(\t))-\pi(q(\t)),
\end{equation}
where $q(\t)$ belongs to $argmax_{q\in Q} \left\{U(\t,q)-\pi(q)\right\}.$ Notice that for all $q$ in $Q$ we have
\begin{equation}\label{Max-ut2}
v(\t)\ge U(\t, q)-\pi(q)
\end{equation}
Analogous to the concepts of {\sc subdifferential}  and {\sc convex conjugate} from classical Convex Analysis, we have that the subset of $Q$ where (\ref{Max-ut2}) is an equality is called the $U$-{\sc subdifferential} of $v$ at $\t$ and $v$ is the $U$-{\sc conjugate} of $
\pi$ (see, for example, \cite{kn:Car}). We write

\[
v(\t)=\pi^U(\t)
\]
and
\[
\partial_Uv(\t):=\{q\in Q\,\mid\,\pi^U(\t)+\pi(q)=U(\t, q)\}
\]
To simplify notation let $\pi(q(\t))=\pi(\t).$ A single pair $(q(\t), \pi(\t))$ is called a {\sc contract}, whereas $\{(q(\t),\pi(\t))\}_{\t\in\T}$ is called a {\sc catalogue}. A catalogue is called {\sc incentive compatible} if $v(\t)\ge v_0(\t)$ for all $\t\in\T,$ where $v_0(\t)$ is type's $\t$ non-participation (or reservation) utility. We normalize the reservation utility of all agents to zero, and assume there is always an {\sc outside option} $q_0$ that denotes non-participation. Therefore we will only consider functions $v\ge 0.$ The Principal's aim is to  devise a pricing function $\pi:Q\to\re$ as to maximize her income

\begin{equation}\label{eq:profit}
\int_{\T}\left(\pi(\t)-C(q(\t))\right)f(\t)d\t
\end{equation}
Inserting (\ref{Max-ut}) into (\ref{eq:profit}) we get the alternate representation

\begin{equation}\label{eq:pprofit}
\int_{\T}\left(U(\t, q(\t))-v(\t)-C(q(\t))\right)f(\t)d\t.
\end{equation}
Expression (\ref{eq:pprofit}) is to be maximized over all pairs $(v, q)$ such that $v$ is U-convex and non-negative and $q(\t)\in\partial_Uv(\t).$ Characterizing $\partial_Uv(\t)$ in a way that makes the problem tractable can be quite challenging. In the case where   $U(\t, q(\t))=\t\cdot q(\t),$ as in \cite{kn:rch}, for a given price schedule $\pi:Q\to\re,$ the maximal net utility
of an agent of type $\t$ is

\begin{equation}\label{eq:v}
v(\t):=\max_{q\in Q}\left\{ \t\cdot q-\pi(q)\right\}
\end{equation}
Since $v$ is defined as the supremium of its affine
minorants, it is a convex function of the types. It follows from the
Envelope Theorem that the maximum  in equation (\ref{eq:v}) is
attained if $q(\t)=\nabla v(\t),$  and we may write

\begin{equation}\label{v-convex}
v(\t)= \t\cdot\nabla v(\t)-\pi(\nabla v(\t)).
\end{equation}
The principal's aggregate surplus is given by

\begin{equation}\label{eq:surplus}
\int_{\T}\left(\pi (q(\t))-C(q(\t))\right)f(\t)d\t.
\end{equation}
Inserting (\ref{v-convex})  into (\ref{eq:surplus}) we get that  the principal's objective is to maximize

\begin{equation}\label{Ppal-Conv}
I[v]:=\int_{\T}\left(\t\cdot\nabla v(\t)-C(\nabla v(\t))-v(\t)\right)f(\t)d\t
\end{equation}
over the set

\[
{\cal{C}}:=\left\{v:\T\to\re\,\mid\, v\,{\mbox{convex}},\, v\ge 0,\,\nabla v(\t)\in Q\right\}.
\]

\subsection{The Musa-Rosen Problem in a Square}
The following structures are shared in the first two examples:

\begin{itemize}

\item $x=(v^k, D^k),$  this structure will determine any possible
 candidate for a minimizer to $J_k(\cdot,\cdot)$
 in the following way:  $v^k$ is a vector of length $k^2$ that will
 contain the approximate values of the optimal function $\v$  evaluated on the points of
 the lattice. The vector $D^k$ has length $2*k^2$ and it contains
 what will be the partial derivatives of $\v$ at the same points

\item $h$ is a vector of length $3*k^2.$ The product $hx$ provides the
discretization of the integral $\int_{\T}(\t\cdot\nabla {v}-{v}(\t))f(\t)d\t.$

\item $B$  is the matrix of constraints. The inequality $Bx\le 0$ imposes the
non-negativity of $v$ and $D$ and the convexity of the resulting $\v_k.$
\end{itemize}

\begin{rem} The density $f(\t)$ is  ``built into" vector $h$ and the cost function $C.$\end{rem}

\begin{figure}[ht!]
\begin{center}
\subfigure[\label{figure:vvv} The optimal function $v$.]
{\epsfig{figure=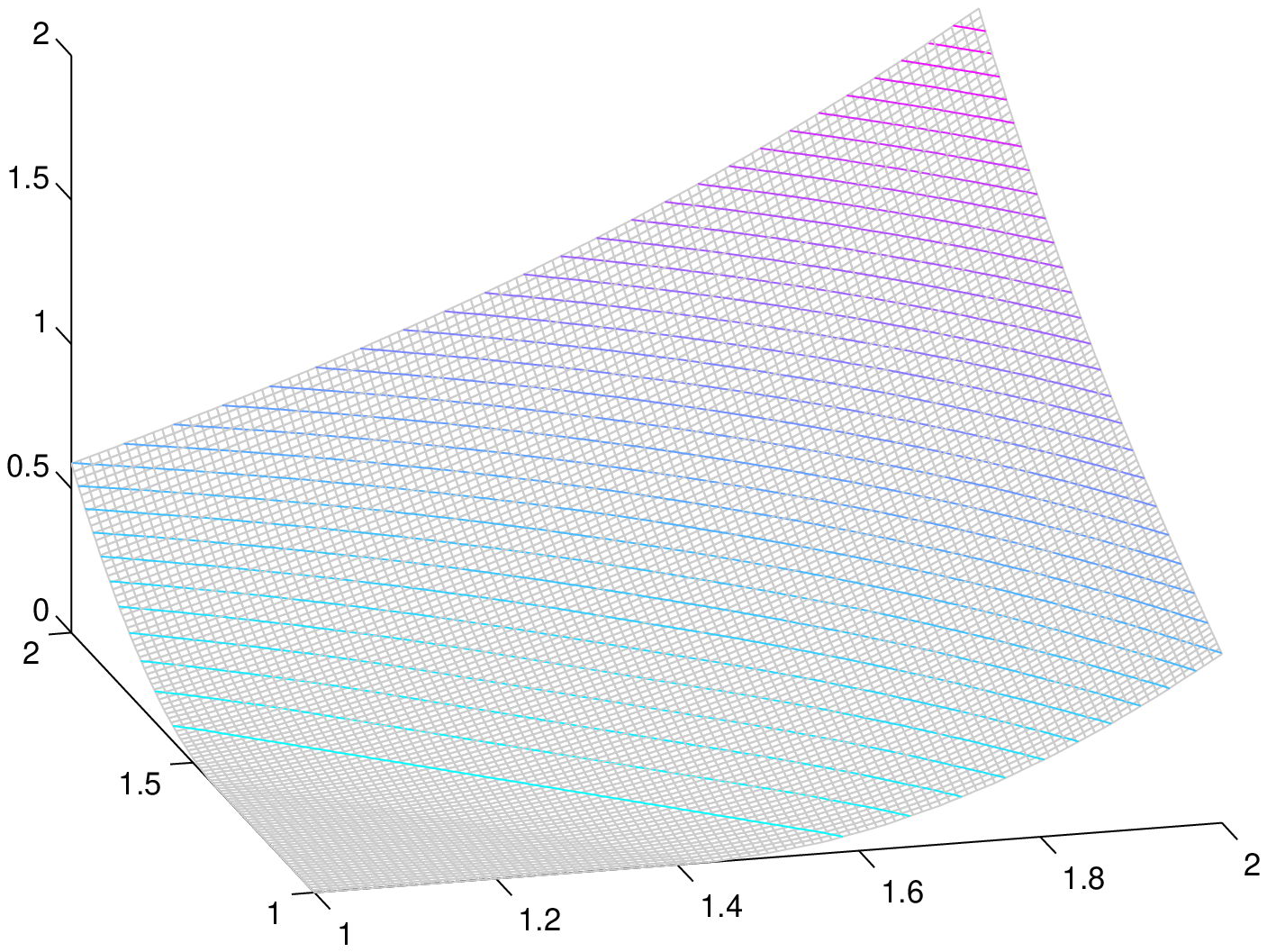,scale=0.45}} \hspace{5mm}
\subfigure[\label{figure:v1}The qualities traded.]
{\epsfig{figure=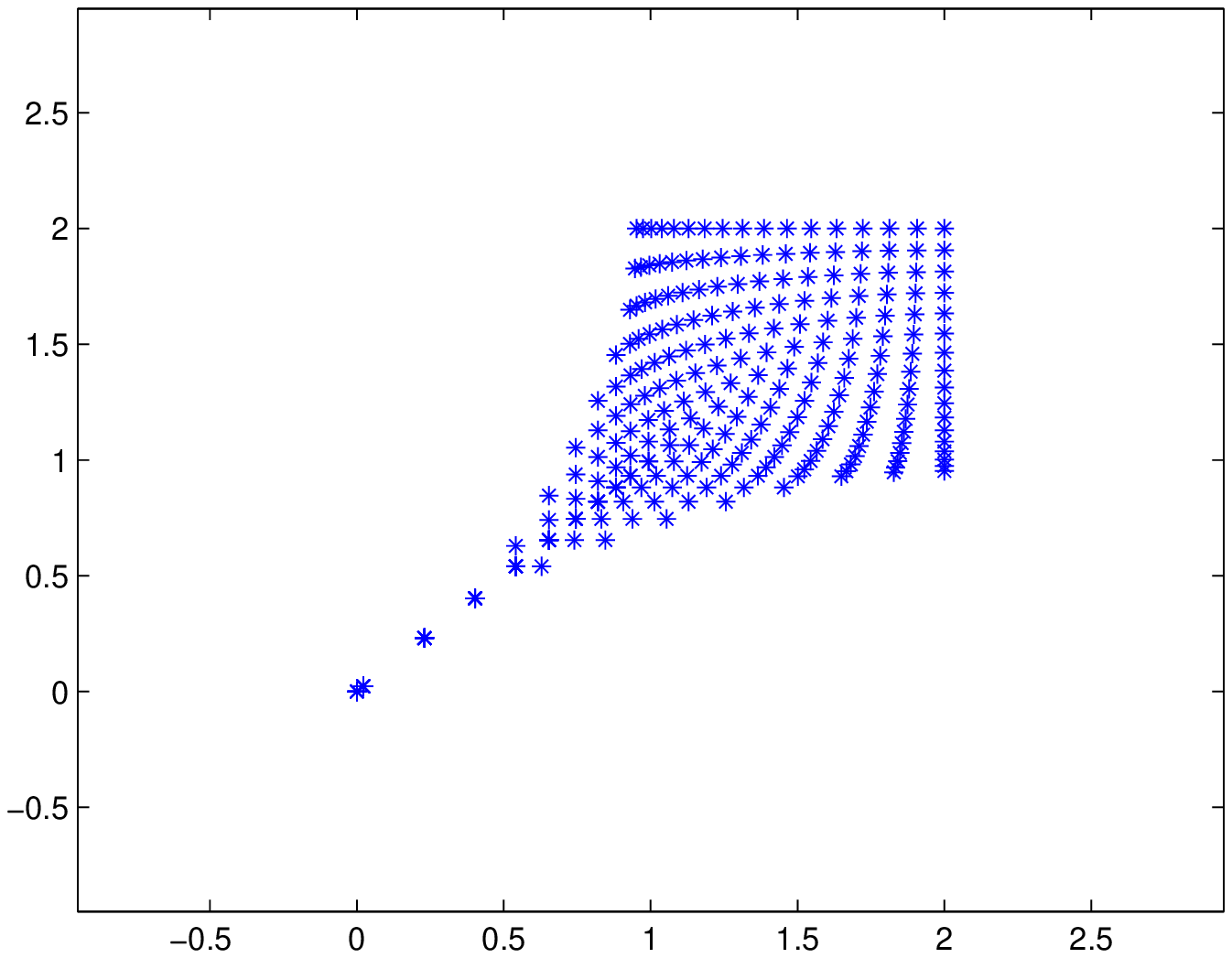,scale=0.45}} \caption{\label{figure:1}
Optimal solution for Uniformly distributed agent types.}
\end{center}
\end{figure}

Let $\T=[1,2]^2,$  $C(q)=\frac{1}{2}\|q\|^2,$ and assume the types
are uniformly distributed. This is our the benchmark problem, since
the solution to the principal's problem  can be found explicitly
\cite{kn:rch}. In this case we have to solve the quadratic program

\[
\sup_x hx-\frac{1}{2}x^tHx
\]
subject to
\[
Bx\le 0
\]
$H$ is a $(3*k^2)\times(3*k^2)$ matrix whose first $k^2$ columns
are zero, since $v$ does not enter the cost function; the four
$k^2\times k^2$ blocks towards its lower right corner form a
$(2*k^2)\times (2*k^2)$ identity matrix. Therefore $\frac{1}{2}x^tHx$ is
a discretization of $\int_{\T}\|\nabla {v}(\t)\|d\t.$ Figure \ref{figure:vvv} was produced using a $17\times 17$-points lattice
and a uniform density, whereas Figure \ref{figure:v1} shows the traded qualities.

\subsection{The Musa-Rosen Problem with a Non-uniform Density}

\begin{figure}[ht!]
\begin{center}
\subfigure[\label{figure:vvvv} The optimal function $v$.]
{\epsfig{figure=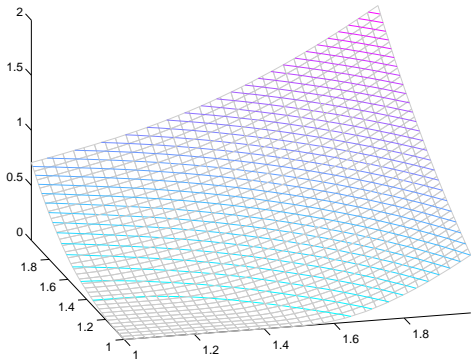,scale=0.45}} \hspace{5mm}
\subfigure[\label{figure:v1}The qualities traded.]
{\epsfig{figure=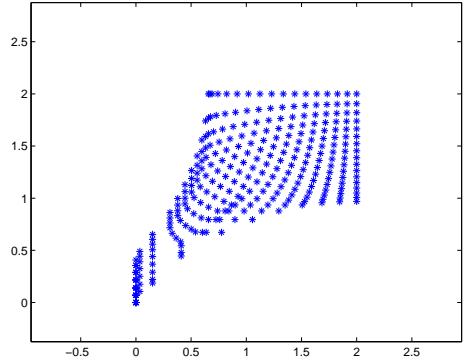,scale=0.45}} \caption{\label{figure:2}
Optimal solution for normally distributed agent types.}
\end{center}
\end{figure}

We keep the cost function of the previous example, but now assume
the types are distributed according to a bivariate normal
distribution with mean $(1.9, 1)$ and variance-covariance matrix

\[
\left[\begin{array}{ll}
.3 & .2\\
.2 & .3
\end{array}\right]
\]
As noted before, the weight assigned to each to each agent
type is built into $h$ and $H,$ so the vector $x$ remains unchanged.
We obtain figure \ref{figure:vvvv}.

\begin{rem} It is interesting to see that in this case bunching of the second kind, as described by Roch\'e and Chon\'e in \cite{kn:rch},
appears to be eliminated as a consequence of the skewed distribution of the agents. This can be seen in the non-linear level curves
of the optimizing function $v.$ This is also quite evident in the plot of the qualities traded, which is shown below.
\end{rem}

The MatLab programs for the two previous examples were run on MatLab 7.0.1.24704 (R14) in a Sun Fire V480 (4$\times$1.2 HGz Ultra III, 16 GB RAM) computer running Solaris 2.10 OS. In the first example 57.7085 seconds of processing time were required. The running time in the second example was 81.7280 seconds.

\subsection{An example with Non-quadratic Cost}

In this example we approximate a solution to the problem of a principal who is selling over-the-counter financial derivatives to a set of heterogeneous agents. This model is presented by Carlier, Ekeland and Touzi in \cite{kn:cet}. They start with a standard probability space $(\O, {\cal F}, P),$ and the types of the agents are given by their risk aversion coefficients under the assumption of mean-variance utilities; namely, the set of agent types is $\T=[0,1],$ and the utility of an agent of type $\t$ when facing product $X$ is

\[
U(\t, X)=E[X]-\t{\mbox Var}[X]
\]

\begin{figure}[ht!]\label{figure:cet}
\begin{center}
\subfigure
{\epsfig{figure=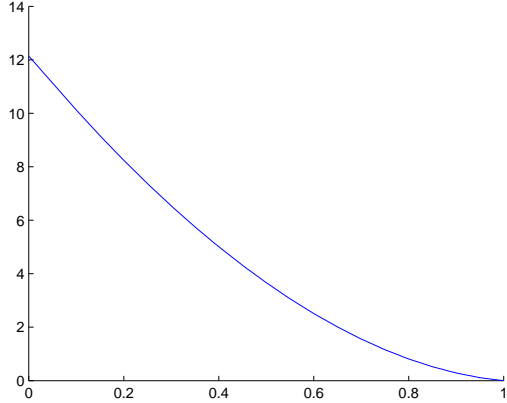,scale=0.55}}  \caption{\label{figure:3}
The optimal function $\v$.}
\end{center}
\end{figure}

Under the assumptions of a zero risk-free rate and that the principal has access to a complete market, her cost of delivering product $X(\t)$ is given by $-\sqrt{-\xi v'(\t)};$ where $\xi$ is the variance of the Radon-Nikodym derivative of the (unique) martingale measure, and Var$[X(\t)]=-v'(\t).$ The principal's problem can be written as

\begin{equation}
\sup_{v\in{\cal C}}\int_{\T}\left(\t v'(\t)+\sqrt{-v'(\t)}-v(\t)\right)d\t
\end{equation}

\noindent where ${\cal C}:=\left\{v:\T\to\re\,\mid\,v\,\,{\mbox {convex}},\,v\ge 0, v'\le 0\,\, {\mbox{ and Var}}\,\,[X(\t)]=-v'(\t)\right\}.$ Figure \ref{figure:cet} shows an approximation of the maximizing $\v$ using 25 agent types.

\subsection{Minimizing Risk}\label{Risk}

The microeconomic motivation for this section is the model of Horst \& Moreno \cite{kn:HM}. We present an overview for completeness.  The principal's income, which is exposed to non-hedgeable risk factors, is represented by $W \leq 0.$ The latter is a bounded random variable defined
on a standard, non-atomic, probability space $\left(\Omega,{\cal{F}}, \Prob \right).$  The principal's goal is to lay off parts of her risk with
the agents whose preferences are mean-variance. The agent types are indexed by their coefficients of
risk aversion, which are assumed to lie $\Theta = [a,1]$ {for some } $a>0.$  The principal underwrites call options on her income with type-dependent strikes:

\[
    X(\t) = (|W|-K(\t))^+ \quad \mbox{with} \quad 0 \leq K(\t)
    \leq \|W\|_\infty.
\]

If the principal issues the catalogue $\{(X(\t), \pi(\t))\}$, she receives a
cash amount of $\int_{\T} \pi \left(\theta \right)d  \theta$
and is subject to the additional liability $\int_{\T} X(\theta)
d\t.$ She evaluates the risk associated with her overall
position
\[
    W + \int_{\T}(\pi(\t)-X(\t))d\t
\]
via the ``entropic measure" of her position, i.e.
\[
\rho\left(W+\int_{\T}(X(\t)-\pi(\t))d\t\right)
\]
where $\rho(X)=log E[exp\{-\beta X\}]$ for some $\beta>0.$  The principal's problem is to devise a catalogue as to minimize her risk exposure. Namely, she chooses a function $v$ and contracts $X$ from the set

\[
    \{(X,v) \mid v \in {\cal C}, \, v \leq K_1, \, -\textnormal{Var}[(|W|-K(\t))^+] =
    v'(\t), \, |v'| \leq K_2, \, 0 \leq K(\t) \leq \|W\|_\infty \},
\]
in order to minimize
\[
    \rho\left( W - \int_{\T}\left\{ (|W|-F(v'(\t)))^+ - \E[(|W|-F(v'(\t)))^+] \right\} d\right) -
    I(v).
\]
where
\[
   I(v) = \int_{\T} \left( \t v'(\t) - v(\t) \right) d\t.
\]
We assume the set of states of the World is finite
with cardinality $m.$  Each possible state $\omega_j$  can occur
with probability $p_j.$ The realizations of the principal's wealth
are denoted by $W=(W_1,\ldots,W_m).$ Note that $p$ and $W$ are
treated as known data. The objective function of our non-linear program is

\begin{eqnarray*}
    F(v,v',K) &=& log\left(exp\left\{-\sum_{i=1}^nW_ip_i+\frac{1}{n}\sum_{i=1}^n\left(\sum_{j=1}^n
    T(K_j-|W_i|)\right)p_i \right.\right.\\
    &-&\left.\left.\frac{1}{n}\sum_{i=1}^n\left(\sum_{j=1}^{n}T(K_j-|W_i|)\right)p_i \right\}\right)\\
    &+&
    \frac{1}{n}\sum v_i-\t_i v_i'
\end{eqnarray*}
where $K=(K_1,\ldots,K_n)$ denotes the vector of type dependent strikes.
We denote by $ng$ the total number of constraints. The principal's problem is to find

\[
    \min_{(v, v', K)} F(v, v', K) \quad \mbox{subject
    to} \quad G(v, v', K)\le 0
\]
where $G:\re^{3n}\to\re^{ng}$ determines the constraints that
keep $(v,v',K)$ within the set of feasible contracts. Let $(1/6, 2/6,\ldots, 1)$ be the uniformly distributed agent types, and

\begin{itemize}

\item $W=4*(-2, -1.7, 1.4,-.7, -.5, 0),$

\item $P=(1/10, 1.5/10, 2.5/10, 2.5/10, 1.5/10, 1/10).$

\end{itemize}
The principal's initial evaluation of her risk is $1.52$. The following are the plots for the approximating $\v$ and the strikes:

\begin{figure}[ht!]
\begin{center}
\subfigure[\label{OptimalV} The optimal function $v$.]
{\epsfig{figure=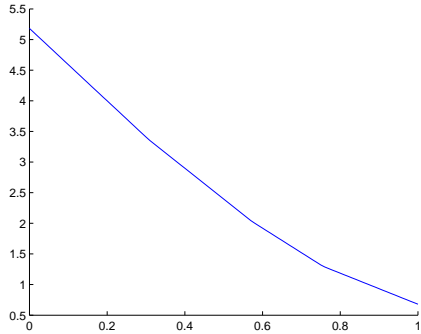,scale=0.45}} \hspace{5mm}
\subfigure[\label{Stikes}The type-dependent strikes ]
{\epsfig{figure=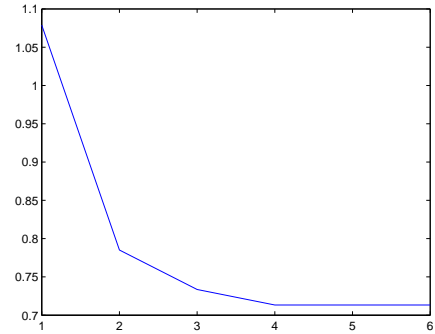,scale=0.45}} \caption{\label{figure:3}
Optimal solution for underwriting call options.}
\end{center}
\end{figure}
Note that for illustration purposes we have changed the scale for the agent types in the second plot. The interpolates of the approximate to the optimal function $\v$ and the strikes are:

\begin{center}
\begin{tabular}{|c|c|}
   \hline
  $\v_{1}$ & 4.196344 \\\hline
  $\v_{2}$ & 3.234565 \\\hline
  $\v_{3}$ & 2.321529 \\\hline
  $\v_{4}$ & 1.523532 \\\hline
  $\v_{5}$ & 0.745045 \\\hline
  $\v_{6}$ & 0.010025\\\hline
\end{tabular}\hspace{1in}
\begin{tabular}{|c|c|}
   \hline
  $K_{1}$ & 1.078869 \\\hline
  $K_{2}$ & 0.785079 \\\hline
  $K_{3}$ & 0.733530 \\\hline
  $K_{4}$ & 0.713309 \\\hline
  $K_{5}$ & 0.713309 \\\hline
  $K_{6}$ & 0.713309 \\\hline
\end{tabular}\end{center}
The Principal's valuation of her risk after the exchanges
with the agents decreases from $11.49$ to  $-3.56.$

\begin{rem}Notice the ''bunching" at the bottom.\end{rem}

\section{Conclusions}
In this paper we have developed a numerical algorithm to estimate the minimizers of variational problems with convexity constraints, with our main motivation stemming from Economics and Finance. Ours is an {\sc internal} method, so at each precision level the approximate minimizers lie within the acceptable set of functions. Our examples are developed over one or two dimensional sets for illustration reasons, but the algorithm can be implemented in higher dimensions. However, it must be mentioned that, as is the case with the other methods found in the related literature, implementing convexity has a high computational cost which increases geometrically with dimension.

\section{Appendix}
\begin{appendix}

\section{Some technical results}\label{ConvergenceConvex}

In order to prove convergence of our algorithm we make use of the Convex Analysis results contained in this section. We will work on $C,$ an open and convex subset of $\re^n.$ 

\begin{definition}A mapping $F:C\to P(\re^m)$ (the power set of $\re^m$)  is said to be {\sc set valued} if for each $x$ in $C,$ $F(x)$ is a non-empty subset of $\re^m.$ 
\end{definition}
Recall that if $f:C\to\re$ is a convex function, then the {\sc subdifferential} of $f$ at $x,$ defined as
\[
\partial f(x):=\{b\in\re^n\,\mid\,f(y)-f(x)\ge b\cdot (y-x)\quad{\mbox{for all}}y\in C\},
\]
is a non-empty subset of $\re^n$ for all $x$ in $C.$ Therefore, the set valued mapping
\[
x\to\partial f(x)
\]
is well defined on $C.$ Notice that if for some $x$ in $C$ we have $\#\partial f(x)=1,$ then
\[
\partial f(x)=\{\nabla f(x)\}.
\]
In such case we say the {\sc subdifferential mapping} is single valued at $x$ and we can simply identify it with the gradient of $f$ at $x.$

\begin{definition} \label{diff}Let $F:C\to\re^m$ be a set valued function.  Then we say $F$ is {\sc differentiable} at $x_0$ iff there is a linear mapping $L_{x_0}:C\to\re^m$ such that for all $\epsilon>0$ there is $\delta>0$ such that if $y_0\in f(x_0),\,y\in f(x)$ and $\|x-x_0\|<\delta$ then
\[
\frac{\|y-y_0-L_{x_0}(x-x_0)\|}{\|x-x_0\|}\le\epsilon
\]
\end{definition}
The following theorem is due to Alexandrov ( \cite{kn:HW})

\begin{thm} Let $f:\re^n\to\re$ be convex, then se set valued function $\partial f$ is differentiable almost everywhere.
\end{thm}
Clearly, in the case where $f$ is single valued,  definition \ref{diff} is equivalent to the regular definition of a Fr\'echet differentiable function. Moreover if $f$ differentiable at $x_0$ and we choose $y_0,$ $y_1$ in $f(x_0)$ and let $x=x_0$ in definition  \ref{diff}  we get
\[
\|y_0-y_1\|\le 0,
\]
which implies $f$ is single valued at $x_0.$ It follows from Alexandrov's Theorem and the observation above that for almost all $\t\in\T,$ the set valued mapping $\t\to\partial f(\t)$ can be identified  with the single valued assignment $\t\to\nabla f(\t),$ and we have the following

\begin{cor}\label{lm:ConvCont} Let $f:C\to\re$ be convex . Then the mapping
\[
\t\to\nabla f(\t)
\]
is well defined and continuous almost everywhere.
\end{cor}

\begin{prop}\label{pr:vprime} Let $A\subset\re^n$ be a convex, open set. Assume the sequence of convex functions $\{f_k:A
to\re\}$ converges uniformly to $\bar{f},$ then $\nabla f_k\to\nabla\bar{f}$
almost everywhere on $A.$
\end{prop}
\begin{Proof}  Denote by $D_if$ the derivative of $f$ in the direction of
$e_i.$ The convexity of $f_k$ and $\bar{f}$ implies the existence of
a set $B,$ with $\mu(A\setminus B)=0$ such  that the partial
derivatives of $f_k$ and $\bar{f}$  exist and are continuous in $B.$
Let $x\in B.$ To prove that $D_i f_k(x)\to D_i \bar{f}(x),$ consider
$\eta$ such that $x+\eta e_i\in A.$ Since $f_k$ is convex

\[
\frac{f_k(x+he_i)-f_k(x)}{h} \ge D_i f_k(x)
\ge\frac{f_k(x-he_i)-f_k(x)}{h}
\]
for all $0<h<\eta.$ Hence
\[
\frac{f_k(x+he_i)-f_k(x)}{h} - D_i\bar{f}(x) \ge D_i f_k(x) -
D_i\bar{f}(x) \ge\frac{f_k(x-he_i)-f_k(x)}{h} - D_i\bar{f}(x).
\]
The left-hand side of this inequality is equal to
\[
\frac{f_k(x+he_i)-\bar{f}(x+he_i)}{h}+\frac{\bar{f}(x)-f_k(x)}{h}+\frac{\bar{f}(x+he_i)-\bar{f}(x)}{h}
- D_i\bar{f}(x).
\]
For $\epsilon >0$ let $0<\delta<\eta$ be such that
\[
\left| \frac{\bar{f}(x+he_i)-\bar{f}(x)}{h} - D_i\bar{f}(x)\right|
<\epsilon
\]
for $|h|\le\delta.$ Let $N\in\n$ be such that

\[
-\epsilon\delta\le f_n(x)-\f(x)\le \epsilon\delta
\]
$n\ge N.$  Hence, taking $h=\delta,$ we have that for all
$n\ge N,$

\[
\dfrac{f_n(x+he_1)-\f(x+he_1)}{h}\le\epsilon\quad{\mbox{and}} \quad\dfrac{\f(x)-f_n(x)}{h}\le\epsilon.
\]
Hence
\[
3\epsilon \ge D_1 f_n(x)-D_1\f(x)
\]
for all $x\in B.$ The same argument shows that
\[
-3\epsilon \le D_1 f_n(x)-D_1\f(x)
\]
for all $n\ge N$ and all $x\in B,$ which concludes the proof.
\end{Proof}
\begin{prop}\label{convex:C1} Let $U\subset\re^n$ be a convex, compact set and let $g:U\to\re$ be a  convex function such that for all $x\in U,$ the subdifferentials $\partial g(x)$ are contained in $Q$ for some compact set $Q.$ Then there exists $\{g_j:U\to\re\}$ such that $g_j\in C^1(U)$ and
$g_j\to g$ uniformly on U.
\end{prop}
\begin{Proof}
Fix $\delta>0$ and define
\[U_{\delta}:=\{(1+\delta)x\,\mid\,x\in U\}.\]
Extend $g$ to be defined on $U_{\delta}.$ Let $K_{\epsilon}$ be a family of mollifiers (see, for instance \cite{kn:GH}), then the functions
\[h_{\epsilon}:=g*K_{\epsilon}\]
are convex, smooth and they converge uniformly to $g$ on $U$ as long as $\epsilon$ is small enough so that
\[U_{\epsilon}:=\{x\in U_{\delta}\,\mid\,d(x,\partial U_{\delta})>\epsilon\}\]
is contained in $U.$ Let $n\in\n$ be such that $U_{1/n}\subset U,$ then the sequence $\{g_j:=h_{1/j}\}$ has the required properties.
\end{Proof}

\begin{lemma}\label{eq:affineconv}
Consider $\phi(\x, z, p)\in C^1\left(\T\times\re\times
Q\to\re\right),$ where $\T=[a,b]^n$ and $Q$ is a compact convex
subset of $\re^n.$ Let $\{f_k:\T\to\re\}$ be a family of convex functions such that $\partial f_K(\t)\subset Q$ for all $\t\in\T,$ and whose uniform limit is $\f.$  Let $\Sigma_k$ be the uniform partition of $\T$ consisting of  $k^n$  cubes of volume $\|\Sigma_k\|:=\left(\frac{b-a}{k}\right)^n.$ Denote by $\sigma_j^k,$ $1\le j\le k^n,$ be the elements of $\Sigma_k$ and let
\[
\|\Sigma_k\|\sum_{i=1}^{k^n} \phi(\x_j^k, f_k(\x_j^k), \nabla f_k(\x_j^k))
\]
be the corresponding Riemann sum approximating $\int_{\T}\phi(x, f_k(\x), \nabla f_k(\x))d\x,$ where $\x_j^k\in\sigma_j^k$ and $\sigma_j^k\in\Sigma_k.$ Then for any $\epsilon>0$ there is $K\in\n$ such that

\begin{equation}\label{eq:2}
\left|\int_{\T}\phi(\x, f_k(\x), \nabla f_k(\x))d\x-\|\Sigma_k\|\sum_{i=1}^{k^n} \phi(\x_j^k, f_k(\x_j^k), \nabla f_k(\x_j^k))\right|\le\epsilon
\end{equation}
for any $k\ge K.$
\end{lemma}
\begin{Proof}  By lemma \ref{convex:C1}, for each $f_k$ there exists a sequence of continuously differentiable convex functions $\{g^k_j\}$ such that
\[
g^k_j\to f_k\quad{\mbox{uniformly}}.
\]
Let $h_k$ be the first element in $\{g^k_j\}$ such that $\|h_k-f_k\|\le\frac{1}{k}$ and $\|\nabla h_k(\t)-\nabla f_k(\t)\|\le\frac{1}{k}$ for all $\t\in\T$ where $\nabla f_k$ is continuous.  Then $h_k\to\f$ uniformly, and by Lemma \ref{pr:vprime} we have that $\nabla h_k(\t)\to\nabla\f(\t)$  a.e. It follows from Egoroff's theorem that for every $n\in\n$ there exists a set $\Lambda_{n}\subset\T$ such that:
\[
\mu(\T\setminus\Lambda_{n})<1/n\quad{\mbox{and}}\quad\nabla h_k\to\nabla\f\quad{\mbox {uniformly on}}\quad\Lambda_{n}.
\]
Let $\X_{\sigma_j^k}(\cdot)$ be the indicator function of $\sigma_j^k$ and define
\[
g_k(\t):=\phi(\t, h_k(\t), \nabla k_k(\t))-\sum_{j=1}^{k^n}\X_{\sigma_j^k}(\t)\phi(\t_j^k, h_k(\t_j^k), \nabla h_k(\t_j^k))
\]
Fix $n,$ consider $\t_0\in\Lambda_n$ and let $\{\t_{0}^k\}$ be the sequence of $\t_j^k$'s converging to $\t_0$ as the partition is refined. By uniform convergence, $\nabla\f$ is continuous on $\Lambda_n,$ hence

\begin{equation}\label{eq:pointwise}
\lim_{k\to\infty} h_k(\t_{0}^k)=\f(\t_{0})\quad{\mbox{and}}\quad\lim_{k\to\infty} \nabla h_k(\t_{0}^k)=\nabla\f(\t_{0})
\end{equation}
It follows from (\ref{eq:pointwise}) and the continuity of $\phi$ that $g_k\to 0$ almost everywhere on $\Lambda_n.$ Notice that as a consequence of the compactness of $\T$ and $Q$ and the definition of $h_k$ we have
\[
\|\phi(\x, f_k(\x), \nabla f_k(\x))\|\le K_1,\quad{\mbox{for al}}\quad\x\in\T\quad{\mbox{and some}}\quad K_1>0.
\]
and
\[
\left| g_k(\t)-\left(\phi(\t, f_k(\t), \nabla f_k(\t))-\sum_{j=1}^{k^n}\X_{\sigma_j^k}(\t)\phi(\t_j^k, f_k(\t_j^k), \nabla f_k(\t_j^k))\right)\right|\le\frac{K_2}{k}
\]
for some $K_2>0$ and all $\t\in\T$ where $\nabla f_k$ is continuous. Therefore
\begin{equation}
\left|\int_{\T}\phi(\x, f_k(\x), \nabla f_k(\x))d\x-\|\sigma_j^k\|\sum_{i=1}^{k^n} \phi(\x_j^k, f_k(\x_j^k), \nabla f_k(\x_j^k))\right|\le\frac{K_2\|\T\|}{k}+\left|\int_{\T}g_k(\t)d\t\right|
\end{equation}
By Lebesgue Dominated Convergence
\[
\lim_{k\to\infty}\int_{\Lambda_{n}}g_k(\t)d\t=0,
\]
moreover, the definition of $\Lambda_n$ implies
\[
\int_{\T\setminus\Lambda_{n}}g_k(\t)d\t\le\frac{2K_1}{n}.
\]
Given $\epsilon>0$ take $n\in\n$ such that $\frac{2K_1}{n}\le\frac{\epsilon}{2}$ and $K$ such that
\[
\frac{K_2\|\T\|}{K}+\left|\int_{\Lambda_{n}}g_K(\t)d\t\right|\le\frac{\epsilon}{2}.
\]
Then equation (\ref{eq:2}) holds for all $k\ge K.$
\end{Proof}

\section{MatLab code for the examples in section \ref{Examples}}

\end{appendix}


\end{document}